\documentclass[12pt]{amsart}
\input colordvi

\usepackage{amssymb}
\usepackage{amsmath}

\newcommand{\Rb}[1]{{\mathbb{R}_{#1}}}
\newcommand{\DEQS}{\begin{eqnarray*}}
\newcommand{\EEQS}{\end{eqnarray*}}
\newcommand{\DEQSZ}{\begin{eqnarray}}
\newcommand{\EEQSZ}{\end{eqnarray}}

\newcommand{\lk}{\left}
\newcommand{\lqq}{\lefteqn}
\newcommand{\rk}{\right}

\newcommand{\CE}{{{ \mathcal E }}}

\newcommand{\CS}{{{ \mathcal S }}}

\newcommand{\CB}{{{ \mathcal B }}}
\newcommand{\EE}{\mathbb{E}}
\newcommand{\PP}{\mathbb{P}}
\newcommand{\CM}{{{ \mathcal M }}}
\newcommand{\CP}{{{ \mathcal P }}}
\newcommand{\CF}{{{ \mathcal F }}}

\newcommand{\RR}{{\mathbb{R}}}

\newcommand{\DD}{{\rm I \kern -0.2em D}}
\newcommand{\dd}{{\rm I \kern -0.2em D}}
\newcommand{\NN}{{\mathbb{N}}}
\newcommand{\bNN}{{\bar{\mathbb{N}}}}

\newcommand\toup{\nearrow}
\newcommand\todown{\searrow}

\newcommand\vt{\vartheta}
\newcommand\vp{\varphi\,}

\newcommand\del[1]{}
\newcommand\think[1]{}
\newcommand\new[1]{}
\newcommand\zus[1]{}

\newcommand\comd[1]{} 
\newcommand\Redd[1]{} 

\newcommand{\eps}{\varepsilon}

\newcommand{\DS}{\displaystyle}

\newcommand\bdm{\begin{displaymath}}
\newcommand\edm{\end{displaymath}}
\newcommand\bea{\begin{eqnarray}}
\newcommand\eea{\end{eqnarray}}

\newtheorem{theorem}{Theorem}[section]

\newtheorem{Notation}{Notation}
\newtheorem{lem}[theorem]{Lemma}

\newtheorem{prop}[theorem]{Proposition}
\newtheorem{coro}[theorem]{Corollary}
\newtheorem{example}[theorem]{Example}

\newtheorem{remark}{Remark}

\begin{document}

\title[\today Maximal regularity ]
{Maximal regularity for stochastic convolutions driven by Levy
noise}

\author[ZB]{Zdzis{\l}aw  Brze{\'z}niak}
\email{zb500@york.ac.uk}
\address{Department of Mathematics, University of York, Heslington, York YO10
5DD, UK}

\author[EH \today]{
Erika Hausenblas}
\email{ erika.hausenblas@sbg.ac.at}
\address{Department of Mathematics, University of Salzburg,
Hellbrunnerstr. 34, 5020 Salzburg, Austria}

\date{\today}

\thanks{The research of the second named author was supported by a    grant P17273 of the Austrian Science Foundation. The research of the first named author was supported  by a grant. He would like to thank the Department of Mathematics, University of Salzburg, for the hospitality.  The research on this paper was initiated during a visit of both authors to the Centro di Ricerca Matematica Ennio de Giorgi in Pisa (Italy), in July 2006. }

\keywords{Stochastic convolution  \and time homogeneous Poisson random measure  \and
 maximal regularity  \and martingale type $p$ Banach spaces}

\subjclass[2000]{}

%
\maketitle

\begin{abstract}
We show that the result from Da Prato and Lunardi  is valid for stochastic convolutions
 driven by L\'evy processes.
\end{abstract}




\section{Introduction}

The aim of the article is to investigate the maximal regularity of
the Ornstein-Uhlenbeck driven by purely discontinuous noise.
In particular, let $(S,\CS)$ be a measurable space,
$E$ be a Banach space of martingale type $p$, $1< p\le 2$, and
$A$ be an infinitesimal generator of an analytic semigroup $(e ^ {-tA})_{0\le t<\infty}$ in $E$.
We consider the following SPDE written in the It\^o-form
\DEQSZ \label{eq00}
\quad\quad\left\{\begin{array}{rcl}d u(t   ) &=& Au(t - ) \; dt
+ \int_{ S}\xi(t;
                  x) \tilde \eta(dx;dt),
\\ u(0) &=& 0  ,
\end{array}
\right.
\EEQSZ
where $\tilde
\eta$ is a $S$-valued time homogeneous compensated Poisson random measure
defined on  a filtered probability space $(\Omega;\CF;(\CF_t)_{0\le t<\infty};\PP)$
with L\'evy measure $\nu$ on $S$, specified later,
and
$\xi:\Omega\times S\to E$ is a predictable process satisfying certain integrability conditions also specified later.
The solution to \eqref{eq00} is given by the so called Ornstein-Uhlenbeck process
$$
u(t) := \int_0 ^ t \int_{S} e ^ {-A(t-r)}  \xi(r,
                  x)\: \tilde \eta(dx;dr),\quad t>0.
$$
Suppose $1\le q\le p$. Our main result will be the following inequality
\begin{equation}
 \boxed{
 \mathbb{E}  \int_0^T \vert
u(t)\vert_{D_A(\theta+\frac1p,q)} ^p \,dt \leq
C \mathbb{E} \int_0^T  \int_S \vert
\xi(t,z)\vert_{D_A(\theta,q)} ^p \,dt ,}
\end{equation}
where $D_A(\theta,p)$, $\theta\in(0,1)$, denotes the real interpolation space of order $\delta$
between $E$ and $D(A)$.
\par
As mentioned in the beginning, if the Ornstein-Uhlenbeck process is driven by a scalar Wiener process,
the question of maximal regularity was answered by Da Prato in \cite{DaP-83} or Da Prato and
Lunardi \cite{DaP-Lun-98}.
We transfer these results to the Ornstein-Uhlenbeck process driven by purely discontinuous noise.

\begin{Notation}
By $\mathbb{N}$ we denote the set of natural numbers, i.e. $\mathbb{N}=\{0,1,2,\cdots\}$ and by $\bar{\mathbb{N}}$ we denote the set $\mathbb{N}\cup\{+\infty\}$. Whenever we speak about $\mathbb{N}$ (or $\bar{\mathbb{N}}$)-valued measurable functions we implicitly assume that that set is equipped with the trivial $\sigma$-field $2^\mathbb{N}$ (or $2^{\bar{\mathbb{N}}}$).
By $\Rb{+}$ we will denote the  interval
$[0,\infty)$. 
If $X$ is a topological space, then by $\mathcal{B}(X)$ we will denote the Borel $\sigma$-field on $X$.
By $\lambda$ we will denote the Lebesgue measure on $(\mathbb{R},\mathcal{B}(\mathbb{R}))$.
For a measurable space $(S,\CS)$ let $M^ +_S$ be the set of all non negative measures on $(S,\CS)$.

\end{Notation}

\section{Main results}
Suppose that $p\in (1,2]$ and that $E$ is a Banach space of martingale type
$p$. Let $(S,\CS)$ be a measurable space and
 $\nu\in M ^ +_S$. Suppose that $\mathfrak{P}=(\Omega,\CF,(\mathcal{F}_t)_{t\geq 0},\PP)$ is a filtered
probability space, $\eta: \CS\times \CB(\RR_+)\to \bNN$ is time homogeneous Poisson random measure
with intensity measure $\nu$ defined over $(\Omega,\CF,\PP)$ and adapted to filtration $(\mathcal{F}_t)_{t\geq 0}$.
We will denote by $\tilde\eta=\eta-\gamma$ the to $\eta$ associated compensated Poisson random measure
where $\gamma$ is given by
$$
\CB(\RR_+)\times \CS\ni (A,I)\mapsto  \gamma(A,I)=\nu(A)\lambda(I)\in\RR _ +.
$$

We denote by $\CP$ the $\sigma$ field on $\Omega\times \RR_+$ generated by
all sets $A\in \CF\hat \times \CB( \RR_+)$, where $A$ is of the form $A=F\times (s,t]$,
with $F\in\CF_s$ and $s,t\in\RR_+$.
If $\xi:\Omega\times \RR_+\to S$ is $\CP$ measurable, $\xi$ is called predictable.

It is then known, see e.g.\ appendix \ref{appB}, that there exists a unique
continuous linear operator associating with  each predictable
process $\xi:\Rb{+}\times S\times \Omega \to E$ {with}
\begin{equation}
\label{cond-2.01}\mathbb{E} \int_0^T\int_S\vert \xi(r,x)\vert^p\,
\nu(dx)\,dr <\infty, \quad T>0,
\end{equation}
an adapted c\'adl\'ag process, denoted by $\int_0^t \int_S \xi(r,x)
\tilde{\eta}(dx,dr)$, $t\geq 0$ such that if $\xi$ satisfies the above condition (\ref{cond-2.01})
and is a
step process with representation
$$
\xi(r) = \sum_{j=1} ^n 1_{(t_{j-1}, t_{j}]}(r) \xi_j,\quad 0\le r,
$$
where $\{t_0=0<t_1<\ldots<t_n<\infty\}$ is a partition of $[0,\infty)$ and
for all $j$,    $\xi_j$ is  an $\CF_{t_{j-1}}$ measurable random variable, then
\begin{equation}
\label{eqn-2.02} \int_0^t \int_S \xi(r,x)
\tilde{\eta}(dx,dr)=\sum_{j=1}^n  \int_S \tilde \xi_j (x) \eta \lk(dx,(t_{j-1}\wedge t, t_{j}\wedge t] \rk).
\end{equation}

The continuity mentioned above means that there exists a constant
$C=C(E)$ independent of $\xi$ such that

\begin{equation}
\label{ineq-2.03} \mathbb{E} \vert \int_0^t \int_S \xi(r,x)
\tilde{\eta}(dx,dr)\vert ^p \leq C\mathbb{E} \int_0^t\int_S\vert
\xi(r,x)\vert^p\, \nu(dx)\,dr, \; t\geq 0.
\end{equation}

One can prove\footnote{The case $q\in (p,\infty)$ is different and
will be discussed later.}, see e.g.\ the proof of {Proposition 3.3 in \cite{Haus-05-exist}, or Theorem 3.1 in \cite{
Brz-03}} for the case $q<p$, and Corollary \ref{cor_maximal} in Appendix B, that for any $q \in [1,p]$ there
exists a constant $C=C_{q}(E)$ such that for each process $\xi$ as
above and for all $t\geq 0$,
\begin{equation}
\label{ineq-2.04} \boxed{\mathbb{E} \vert \int_0^t \int_S \xi(r,x)
\tilde{\eta}(dx,dr)\vert ^q \leq C\del{C_q(E)}\mathbb{E}
\big(\int_0^t\int_S\vert \xi(r,x)\vert^p\,
\nu(dx)\,dr\big)^{q/p}.}
\end{equation}
\begin{remark}\label{rem-general} Let us denote
\begin{eqnarray*}
I(\xi)(t) &:=& \int_0^t \int_S \xi(r,x) \tilde{\eta}(dx,dr), \;
t\geq
0\\
\Vert \xi\Vert &:=& \left(\int_S\vert \xi(x)\vert^p\,
\nu(dx)\right)^{1/p},\; \xi\in L^p(S,\nu;E).
\end{eqnarray*}

Then the inequality (\ref{ineq-2.04}) takes the following
form
\begin{equation*}
\mathbb{E} \vert I(\xi)(t)\vert ^q \leq C_q(E)\mathbb{E} \left[
\big(\int_0^t\Vert \xi(r)\Vert^p\,dr\big)^{q/p}\right].
\end{equation*}
This should be (and will be) compared with the Gaussian case. Note
that in this case $\Vert \xi\Vert$ is simply the $L^p(S,\nu,E)$ norm
of $\xi$. In the Gaussian case the situation is different.\\
Let us also point out that the inequality (\ref{ineq-2.04}) for
$q<p$ follows from the same inequality for $q=p$. In fact, using
Proposition IV.4.7 from
 \cite{Revuz-Yor-99}, see  the proof of Theorem 3.1 in \cite{Brz-03},
 one can prove a stronger result. Namely that if inequality (\ref{ineq-2.04})
 holds true for $q=p$, then for $q\in [1,p)$  there exists a constant $K_q>0$  such that for
each accessible stopping time $\tau>0$,
\begin{equation}\label{ineq-Burkholder}
\mathbb{E}\sup_{0\le t\le \tau}\vert  I(\xi)(t)\vert ^q\le  K_q\,
\mathbb{E}\,\Bigl( \int^\tau_0\Vert \xi (t)\Vert ^p\, d
t\Bigr)^{q/p}.
\end{equation}
\end{remark}

Assume further that $-A$ is an infinitesimal generator of an
analytic semigroup denoted by $(e^{-tA})_{t\geq 0}$ on $E$.

Define the stochastic convolution of the semigroup $(e^{-tA})_{t\geq
0}$ and an $E$-valued process $\xi$ as above by the following
formula

\begin{equation}
\label{eqn-2.05} SC(\xi)(t)= \int_0^t \int_S e^{(t-r)A}\xi(r,x)
\tilde{\eta}(dx,dr), \; t\geq 0.
\end{equation}

 Let us recall that there exist constants $C_0$ and
$\omega_0$ such that
$$\Vert e^{-tA}\Vert \leq C_0e^{t\omega_0},\; t\geq 0.$$
Without loss of generality, we will assume from now on that
$\omega_0<0$. Let us also recall the following characterization of
the real interpolation\footnote{In order to fix the notation let me
point out that the interpolation functor $(X_0,X_1)_{\theta,q}$,
$\theta \in (0,1)$, $q\in [1,\infty]$, between two Banach spaces
$X_1$ and $X_0$ such that both are continuously embedded into a
common topological Hausdorff vector space, satisfies the following
properties: (i)$(X_1,X_0)_{\theta,q}=(X_0,X_1)_{1-\theta,q}$, (ii)
if $X_0\subset X_1$, $0<\theta_1< \theta_2<1$ and $p,q\in
[1,\infty]$, then $(X_0,X_1)_{\theta_1,p}\subset
(X_0,X_1)_{\theta_2,q}$. Roughly speaking, (ii) implies that, if
$X_0\subset X_1$, then $(X_0,X_1)_{\vartheta,p} \todown X_0$ as
$\vartheta \todown 0$ and $(X_0,X_1)_{\theta,p} \toup X_1$ as
$\vartheta \toup 0$. Or equivalently, if $X_0\subset X_1$, then
$(X_1,X_0)_{\theta,p} \todown X_0$ as $\theta \toup 1$ and
$(X_1,X_0)_{\theta,p} \toup X_1$ as $\theta \todown 1$. See
Proposition 1.1.4 in \cite{Lunardi-99} and section 1.3.3 in
\cite{Triebel-78}.} spaces
$(E,D(A^m))_{\theta,q}=(D(A^m),E)_{1-\theta,q}$, where
$m\in\mathbb{N}$,  between $D(A^m)$ and $E$ with parameters $\theta
\in (0,1)$ and $q\in [1,\infty )$, see section 1.14.5 in
\cite{Triebel-78} or \cite{DaP-83}. If $\delta \in (0,\infty]$ then

\begin{equation*}
(D(A^m),E)_{1-\theta,q} =\left\{ x\in E :\int^{\delta }_{0} |
t^{m(1-\theta) }A^me^{-tA}x| ^{q} {dt\over t} < \infty   \right\}.
\end{equation*}
\begin{equation}\label{eqn-char-int-2.06}
(E,D(A^m))_{\vartheta,q} =\left\{ x\in E :\int^{\delta }_{0} |
t^{m(1-\vartheta) }A^me^{-tA}x| ^{q} {dt\over t} < \infty
\right\}.
\end{equation}

The norms defined by the equality (\ref{eqn-char-int-2.06}) for
different values of $\delta$ are equivalent.

The space $(D(A^m),E)_{1-\theta,q}=(E,D(A^m))_{\theta,q}$ is often denoted by
$D_{A^m}(\theta ,p)$ and we will use the following notation
\begin{equation}\label{eqn-norm-int-2.07}
\boxed{|x|^q_{D_{A^m}(\theta ,q);\delta} =\int^{\delta }_{0} |
t^{m(1-\theta) }A^me^{-tA}x| ^{q} {dt\over t}.}
\end{equation}

In the general case, one has the following equality but only for
$\delta\in (0,\infty)$:

\begin{equation}\label{eqn-char-int-2.08}
(E,D(A^m))_{\theta,q} =\left\{ x\in E :\int^{\delta }_{0} |
t^{m(1-\theta) }(\omega_0I+A)^me^{-t(\omega_0+A)}x| ^{p} {dt\over t}
< \infty  \right\}. \end{equation}

In this case, the formula (\ref{eqn-char-int-2.08}) takes the
following form
\begin{equation}\label{eqn-norm-int-gen-2.09}
|x|^q_{D_{A^m}(\theta ,q);\delta} =\int^{\delta }_{0} |
t^{m(1-\theta) }A^me^{-tA}x| ^{q} {dt\over t}+|x|^q.
\end{equation}
\del{where $\vert\cdot\vert$ denotes the norm in $E$.}

Let us finally recall that if $0<k<m\in\mathbb{N}$, $p\in
[1,\infty]$ and $\theta \in (0,1)$, then
$(E,D(A^k))_{\theta,p}=(E,D(A^m))_{\frac{k}{m}\theta,p}$, see
\cite{Triebel-78} Theorem 1.15.2 (f). Therefore, if $p\in
[1,\infty)$ and $\theta \in[0,1-\frac1p)$, then
\begin{equation}
\label{eqn-reiteration}
D_A(\theta+\frac1p,q)=D_{A^2}(\frac\theta2+\frac1{2p},q)
\end{equation}
with equivalent norms.

 Our main result in this note is the following

\begin{theorem}\label{th-main}
Under the above assumptions, for all \del{$q\in [1,p]$ and} $\theta
\in (0,1-\frac1p)$, there exists a constant
$C=\hat{C}_{\del{q,}\theta}(E)$ such that for any process $\xi$
described above and all $T\geq 0$, the following inequality holds
\begin{equation}
\label{ineq-2.10}\boxed{ \mathbb{E}  \int_0^T \vert
SC(\xi)(t)\vert_{D_A(\theta+\frac1p,q)} ^p\del{q} \,dt \leq
C\del{\hat{C}_{\del{q,}\theta}(E)}\mathbb{E} \int_0^T
\del{\big(}\int_S\vert
\xi(t,z)\vert_{D_A(\theta,q)}^p\,\nu(dx)\del{\big)^{q/p}} \,dt .}
\end{equation}
\end{theorem}

In the Gaussian case and $q=p=2$, and  $E$ being a Hilbert space,
the above result was proved by Da Prato in \cite{DaP-83}. This
result was then generalized  to a class of so called
Banach spaces of martingale type $2$ in \cite{Brz-91-M_type_2}, see also
\cite{Brz-95-M_type_2}, for nuclear Wiener process and in
\cite{Brz-Gat-99}, to the case of
cylindrical Wiener process. Finally, Da Prato and Lunardi studied in
\cite{DaP-Lun-98} the  case when $p=2$ and $q \geq 2$ for a one
dimensional Wiener process. However, a generalisation of the last
result to a cylindrical Wiener process does not cause any serious
problems. We will state corresponding result at the end of this
Note.

Theorem  \ref{th-main} will be deduced from a more general result
whose idea can be traced back to Remark \ref{rem-general}.

\begin{theorem}\label{th-general}  Let $(\Omega,\CF,(\mathcal{F}_t)_{t\geq 0},\mathbb{P})$ is a
filtered probability space, $p\in (1,2]$ and $q\in [p,\infty)$. Let $\mathcal{E}_p$ be a
class of separable Banach spaces satisfying the following
properties.

\begin{enumerate}
\item[(R1)] With each space  $E$ belonging  to the class $\mathcal{E}_p$
we associate  a separable  Banach space $R=R(E)$
such that there is
 a family $(I_t)_{t\geq 0}$ of linear operators from
the  class $\mathcal{M}_{\text{loc}}^p(R(E))$ of all predictable
$R(E)$-valued processes   to
$L^p(\Omega,\mathcal{F}_t,\mathbb{P};E)$ such that
 for some constant $C=C_p>0$
\begin{equation}
\label{ineq-2.11} 
\mathbb{E} \vert I_t(\xi)\vert_E ^p \leq
C_p\mathbb{E}  \big(\int_0^t\Vert \xi(r)\Vert_{R(E)}^p\,dr\big).
\end{equation}
\item[(R2)] If $E \in \mathcal{E}_p$ and $E_1$ isomorphic to $E$,
then $E_1$ belongs to $ \mathcal{E}_p$ as well.
\item[(R3)] If $E_1,E_2 \in \mathcal{E}_p$ and $\Phi:E_1\to E_2$ is a
bounded linear operator, then
$$\Vert \Phi \xi\Vert_{R(E_2)} \leq \vert \Phi\vert \Vert
\xi\Vert_{R(E_1)},\; \xi \in R(E_1).$$
\item[(R4)] If $(E_0,E_1)$ is an interpolation couple such that $E_1,E_2 \in
\mathcal{E}_p$, then the real interpolation spaces
$(E_0,E_1)_{\theta,p}$, $\theta \in (0,1)$, belongs to
$\mathcal{E}_p$ as well.
\item[(R5)] For every $\delta >0$ here exists a constant $K_\delta>0$ such that
\begin{equation}\label{ineq-inter}
\int_0^{\delta}
 \Vert r^{1-\theta} Ae^{-r
 A} \xi \Vert_{R(E)}^p
  \, \frac{dr}{r} \leq K_\delta^p \Vert \xi
  \Vert_{R(D_A(\theta,p))}^p,\; \xi \in
R(E).
\end{equation}
\item[(R6)] There exists
a constant $\hat{C}_q>0$ such that for all $t>0$
\begin{equation}
\label{ineq-2.14} \mathbb{E} \vert I_t(\xi)\vert_E ^q \leq
\hat{C}_q\mathbb{E}  \big(\int_0^t\Vert
\xi(r)\Vert_{R(E)}^p\,dr\big)^{q/p},\; \xi \in
\mathcal{M}_{\text{loc}}^p(R(E)).
\end{equation}
\end{enumerate}
Define another family $(SC_t)_{t\geq 0}$ of linear operators from
$\mathcal{M}_{\text{loc}}^p(R(E))$ to
$L^p(\Omega,\mathcal{F}_t,\mathbb{P};E)$ by the following formula

\begin{equation}
\label{eqn-2.12} SC_t(\xi)= I_t\big( e^{-(t-\cdot)A}\xi(\cdot)
\big), \; t\geq 0.
\end{equation}
Then,  for \del{every $q\in [1,p]$ and} every $\theta \in
(0,1-\frac1p)$, there exists a constant $\hat{C}_{q,\theta}(E)$ such
that for all $T>0$ the following inequality holds
\begin{equation}
\label{ineq-2.13}  \mathbb{E} \int_0^T  \vert
SC_t(\xi)\vert_{D_A(\theta+\frac1p,q)} ^q \,dt \leq
\hat{C}_{q,\theta}(E)\mathbb{E} \int_0^T \Vert
\xi(s)\Vert_{R(D_A(\theta+\frac1p,q))}^q \,dt .
\end{equation}

\end{theorem}
\begin{remark}
It follows from (i) that if $\xi(r)=\eta(r)$ a.s.\ for a.a.\ $r\in
[0,t]$, then $I_t(\xi)=I_t(\eta)$.
\end{remark}

Now we shall present two basic examples.
\begin{example}\label{example-gaussian}
Let $(\Omega,\CF,(\mathcal{F}_t)_{t\geq 0},\mathbb{P})$ be a filtered
probability space, $p=2$. Let $H$ be a separable Hilbert space and
let $\mathcal{E}_2$ be a class of all 2-smoothable Banach spaces.
With $E\in \mathcal{E}_2$ we associate the space $R(E):=R(H,E)$ of
all $\gamma$-radonifying operators from $H$ to $E$. It is known, see
\cite{Neidhardt-78} that $R(H,E)$ is a separable Banach space
equipped with any of the following equivalent
norms\footnote{Equivalence of the norms is a consequence of
Khinchin-Kahane inequality.}, $2\le q<\infty$,
\begin{eqnarray}
\label{eqn-2.14} \Vert \vp   \Vert_{R(H,E);q}^q &:=& \mathbb{E}\vert
\sum_j \beta_j \vp e_j\vert_E^q, \;\vp \in R(H,E),
\end{eqnarray}
 $\{e_k\}_k$ be an ONB of $H$ and $\{\beta_k\}_k$ a sequence of
i.i.d.\ Gaussian N(0,1) random variables.\\
\del{
\textbf{Question} It is of interest to see whether the following is
true. Suppose that $T$ is a linear, densely defined operator from
$H$ to $E$ such that the RHS of (\ref{eqn-2.14}) is finite. Does
this imply that $T$ is a bounded operator from $H$ to $E$?
}
\end{example}

\begin{example}\label{example-poisson}
Let $(\Omega,\CF,(\mathcal{F}_t)_{t\geq 0},\mathbb{P})$ be a filtered
probability space, $p\in(1,2]$.
Let $(S,\CS)$ be a measurable space and $\eta:\CS\hat \times \CB(\RR_+)\to\NN ^ +$ be a time homogeneous, compensated Poisson random measure over $(\Omega;\CF;\PP)$ adapted to filtration $(\CF_t)_{t\ge 0}$
with intensity $\nu\in M^ +_S$.
Let $\CE_p$ be the set of all separable Banach spaces of martingale type $p$. With
$E\in \CE_p$ we associated a measurable transformation $\xi:S\to E$ such that
$$
\int_S |\xi(x)|_E ^ p \nu(dx)<\infty.
$$
Then for $q\in [p,\infty)$ let
$$
\|\xi\| ^ q_{R(E)} := \EE \lk| \int_0 ^ 1 \int_S  \xi(x)\: \tilde \eta(dx,dr)\rk|_E ^ q.
$$
\end{example}

\section{Proof of Theorem \ref{th-general}}
\label{sec:proof_general}

We begin with the case $q=p$. Without loss of generality the norm $\vert \cdot \vert_{D_A(\theta+\frac1p,p);1}$, defined by formula (\ref{eqn-norm-int-2.07}), will be denoted   by
$\vert \cdot \vert_{D_A(\theta+\frac1p,p)}$. Also, we may assume
that $A^{-1}$ exists and is bounded so that the graph norm in $D(A)$
is equivalent to the norm $\vert A\cdot\vert$.

By the equality (\ref{eqn-reiteration}), definition (\ref{eqn-norm-int-2.07}), the  Fubini Theorem and formula
  (\ref{eqn-2.12}) we have
\DEQS
\lqq{ \mathbb{E} \int_0^T  \vert SC_t(\xi)\vert_{D_A(\theta+\frac1p,p)} ^p
\,dt \leq C \mathbb{E} \int_0^T  \vert
SC_t(\xi)\vert_{D_{A^2}(\frac\theta2+\frac1{2p},p)} ^p \,dt
}
&&
\\
&= & C\int_0^T  \int_0^1
\mathbb{E} \vert
r^{2(1-\frac\theta2-\frac1{2p})} A^2e^{-rA} SC_t(\xi)\vert^p\frac{dr}{r}\,dt
\\
&=&C \int_0^T  \int_0^1 r^{p(2-\theta)-1}
\mathbb{E} \vert A^2e^{-{r}A}  I_t\big( e^{-(t-\cdot)A}\xi(\cdot)
\big) \vert^p\frac{dr}{r}\,dt
\\
&=& C\int_0^T  \int_0^1 r^{p(2-\theta)-1}
  \,\mathbb{E} \vert I_t\big(A^2e^{-{r}A} e^{-(t-\cdot)A}\xi(\cdot)
\big) \vert^p\frac{dr}{r}\,dt\leq \cdots
\end{eqnarray*}

By applying next the
inequality (\ref{ineq-2.11}), the property (R3),   the Fubini
Theorem,  the fact that $|Ae^{-\frac{r}2A}|\leq Cr^{-1}$, $r>0$, for
some constant $C>0$ as well as by observing that $\frac1{t-u+r}\leq
\frac1r$ for $t\in [u,T]$, $r>0$, we infer that

\begin{eqnarray*}
\lqq{ \cdots \leq C_p \int_0^1 r^{p(2-\theta)-1}  \int_0^T \mathbb{E}
\int_0^t \Vert   A^2 e^{-(t-u+{r})A}\xi(u) \Vert_{R(E)}^p\,
du \,dt \,\frac{dr}{r}}
\\
&\leq& C_p \int_0^1 r^{p(2-\theta)-1}
\\&& \int_0^T \mathbb{E} \int_0^t
\vert A e^{-\frac{t-u+r}{2}A} \vert^p\,  \Vert A
e^{-\frac{t-u+r}{2}A}\xi(u) \Vert_{R(E)}^p \, du \,dt \,\frac{dr}{r}
\\
&\leq& C_p \mathbb{E}  \int_0^1 r^{p(2-\theta)-1}\lk[
\sup_{0\le u\le t}(t-u+r) ^ {-p} \rk]
\\
&& \int_0^T   \big[
\int_u^{T}  \Vert   A e^{-\frac{t-u+r}{2}A} \xi(u) \Vert_{R(E)}^p\,
dt \big]  \, du \, \frac{dr}{r}
\\
&\leq & C_p \mathbb{E} \int_0^T
\\ && { \int_0^{T+1-\rho}
\Vert
Ae^{-\frac\sigma2
 A} \xi(\rho) \Vert_{R(E)}^p\, \big[ \int_{\rho \vee (\sigma+\rho-1)}^{\rho+\sigma} (\sigma+\rho-\tau)^{p(1-\theta)-2}\,
d\tau\big]
  \, {d\sigma} \, d\rho}
\\
&\leq & C_p \mathbb{E} \int_0^T
\\
&&
\int_0^{T+1-\rho} \Vert
Ae^{-\frac\sigma2
 A} \xi(\rho) \Vert_{R(E)}^p\, \big[ \int_{\rho }^{\rho+\sigma} (\sigma+\rho-\tau)^{p(1-\theta)-2}\,
d\tau\big]
  \, {d\sigma} \, d\rho
\EEQS
\DEQS
&= & C_p \mathbb{E} \int_0^T \int_0^{T+1-\rho} \Vert
Ae^{-\frac\sigma2
 A} \xi(\rho) \Vert_{R(E)}^p\, \big[ \int_{0 }^{\sigma} \tau^{p(1-\theta)-2}\,
d\tau\big]
  \, {d\sigma} \, d\rho
  \\
&= & C_p^\prime \mathbb{E} \int_0^T \int_0^{T+1-\rho}
\sigma^{p(1-\theta)-1} \Vert Ae^{-\frac\sigma2
 A} \xi(\rho) \Vert_{R(E)}^p
  \, {d\sigma} \, d\rho
\\
&\leq & C_p^{\prime\prime} \mathbb{E} \int_0^T \int_0^{T/2}
 \Vert \sigma^{1-\theta} Ae^{-\sigma
 A} \xi(\rho) \Vert_{R(E)}^p
  \, \frac{d\sigma}{\sigma} \, d\rho
  \\ &\leq& \hat{C}_{p}^{\prime\prime\prime} K_T^p\mathbb{E} \int_0^T
 \Vert   \xi(r)
\Vert_{R(D_A(\theta,p))}^p \,  dr,
\end{eqnarray*}
where the last inequality is a consequence of the assumption (R5).
\del{ Applying next inequality
(\ref{ineq-inter}) we infer that
}

The  proof in the case $q>p$ follows the same ideas. Note also that the above
prove resembles closely the proof from \cite{DaP-Lun-98}. We give
full details below.

We consider now the case $q> p$. We use the same notation as in the
previous  case. But we will make some (or the same) additional assumptions.
\del{By the Fubini Theorem and formulae (\ref{eqn-reiteration}),
(\ref{eqn-norm-int-2.07}) and  (\ref{eqn-2.12}) we have}
By the equality (\ref{eqn-reiteration}), definition (\ref{eqn-norm-int-2.07}), the  Fubini Theorem and formula
  (\ref{eqn-2.12}) we have
\begin{eqnarray*}
\lqq{ \mathbb{E} \int_0^T  \vert SC_t(\xi)\vert_{D_A(\theta+\frac1p,q)} ^q
\,dt \leq C \mathbb{E} \int_0^T  \vert
SC_t(\xi)\vert_{D_{A^2}(\frac\theta2+\frac1{2p},q)} ^q \,dt} &&
\\
&=& C\int_0^T  \int_0^1 s^{q(2-\theta)-\frac qp}
\mathbb{E} \vert
A^2e^{-sA} SC_t(\xi)\vert^q\frac{ds}{s}\,dt
\\
&=&C \int_0^T  \int_0^1
s^{q(2-\theta)-\frac{q}{p}} \mathbb{E} \vert A^2e^{-{s}A}  I_t\big(
e^{-(t-\cdot)A}\xi(\cdot) \big) \vert^q\frac{ds}{s}\,dt
\\
&=& C\int_0^T  \int_0^1
s^{q(2-\theta)-\frac{q}{p}}
  \,\mathbb{E} \vert I_t\big(A^2e^{-{s}A} e^{-(t-\cdot)A}\xi(\cdot)
\big) \vert^q\frac{ds}{s}\,dt \leq \cdots
\end{eqnarray*}

Before we continue, we formulate the following simple Lemma.
\begin{lem}\label{lem-simple}
There exists a constant $C>0$ such that for all $t>0$, $s\in (0,1)$
$$\big(\int_0^t \frac1{(t-r+s)^{\frac{pq}{q-p}}}\,dr\big)^{\frac{q}p-1} \leq C \frac1{s^{q(1-\frac1p)+1}}  $$
\end{lem}
\begin{proof}[Proof of Lemma \ref{lem-simple}.]
Denote $\alpha= \frac{pq}{q-p}$ and observe that $\alpha >1$. Since
$\int_0^t \frac1{(t-r+s)^\alpha}\, dr = \int_0^t
\frac1{(r+s)^\alpha}\, dr \leq \int_0^\infty \frac1{(r+s)^\alpha}\,
dr= \frac1{\alpha-1} \frac1{s^{\alpha-1}}$ and
$(\alpha-1)(\frac{q}p-1)=q(1-\frac1p)+1$,  the result follows.
\end{proof}
  As in the earlier case, by  applying  the
inequality (\ref{ineq-2.11}),  the property (R3), the Fubini
Theorem,  the fact that $|Ae^{-\frac{s}2A}|\leq Cs^{-1}$, $s>0$, for
some  constant $C>0$ as well as H\"older inequality and Lemma
\ref{lem-simple} we infer that

\begin{eqnarray*}
\cdots &\leq& \hat{C}_{q} \int_0^1 s^{q(2-\theta)-\frac{q}{p}}
\int_0^T \mathbb{E}\big[ \int_0^t \Vert   A^2 e^{-(t-r+{s})A}\xi(r)
\Vert_{R(E)}^p\,
dr \big]^{q/p}\,dt \,\frac{ds}{s}\\
&\leq& C\hat{C}_{q} \mathbb{E} \int_0^1 s^{q(2-\theta)-\frac{q}{p}}
\int_0^T \big[ \int_0^t \vert A e^{-\frac{t-r+s}{2}A} \vert^p\,
\Vert A e^{-\frac{t-r+s}{2}A}\xi(r) \Vert_{R(E)}^p\, dr \big]^{q/p}
\, dt \, \frac{ds}{s}
\\
&\leq& C\hat{C}_{q} \mathbb{E} \int_0^1 s^{q(2-\theta)-\frac{q}{p}}
\int_0^T \left[ \big(\int_0^t \vert A e^{-\frac{t-r+s}{2}A}
\vert^{\frac{pq}{q-p}}\,dr\big)^{\frac{q}p-1}\,\right.
\\&&\hspace{5truecm}\lefteqn{ \left.\int_0^t \Vert A
e^{-\frac{t-r+s}{2}A}\xi(r) \Vert_{R(E)}^q\, dr \right] \, dt \,
\frac{ds}{s}}
\end{eqnarray*}
\begin{eqnarray*}
&\leq& C^\prime\hat{C}_{q} \mathbb{E} \int_0^1
s^{q(2-\theta)-\frac{q}{p}} \int_0^T \frac1{s^{q(1-\frac1p)+1}}
\,\int_0^t \Vert A e^{-\frac{t-r+s}{2}A}\xi(r) \Vert_{R(E)}^q\, dr
\, dt \, \frac{ds}{s}
\\
&=& C^\prime\hat{C}_{q} \mathbb{E}  \int_0^1 s^{q(1-\theta)-1}
\int_0^T \big[ \int_r^{T} \Vert   A e^{-\frac{t-r+s}{2}A} \xi(r)
\Vert_{R(E)}^q\, dt \big]  \, dr \, \frac{ds}{s}
\\
&\leq & C^\prime\hat{C}_{q} \mathbb{E} \int_0^T \int_0^{T+1-\rho}
\Vert Ae^{-\frac\sigma2
 A} \xi(\rho) \Vert_{R(E)}^q\, \big[ \int_{\rho \vee (\sigma+\rho-1)}^{\rho+\sigma}
 (\sigma+\rho-\tau)^{q(1-\theta)-2}\, d\tau\big]
  \, {d\sigma} \, d\rho
\\
&\leq & C^\prime\hat{C}_{q} \mathbb{E} \int_0^T \int_0^{T+1-\rho}
\Vert Ae^{-\frac\sigma2
 A} \xi(\rho) \Vert_{R(E)}^q\, \big[ \int_{\rho }^{\rho+\sigma} (\sigma+\rho-\tau)^{q(1-\theta)-2}\,
d\tau\big]
  \, {d\sigma} \, d\rho
\\
&= & C^\prime\hat{C}_{q} \mathbb{E} \int_0^T \int_0^{T+1-\rho} \Vert
Ae^{-\frac\sigma2
 A} \xi(\rho) \Vert_{R(E)}^q\, \big[ \int_{0 }^{\sigma} \tau^{q(1-\theta)-2}\,
d\tau\big]
  \, {d\sigma} \, d\rho
  \\
&= & \hat{C}_{q}^\prime \mathbb{E} \int_0^T \int_0^{T+1-\rho}
\sigma^{q(1-\theta)-1} \Vert Ae^{-\frac\sigma2
 A} \xi(\rho) \Vert_{R(E)}^q
  \, {d\sigma} \, d\rho
\\
&\leq & \hat{C}_{q}^{\prime\prime} \mathbb{E} \int_0^T \int_0^{T/2}
 \Vert \sigma^{1-\theta} Ae^{-\sigma
 A} \xi(\rho) \Vert_{R(E)}^q
  \, \frac{d\sigma}{\sigma} \, d\rho \leq \\
&\leq& \hat{C}_{q}^{\prime\prime}K_{T/2}^p \mathbb{E} \int_0^T
 \Vert   \xi(r)
\Vert_{R(D_A(\theta,q))}^q \,  dr
\end{eqnarray*}
where  the last inequality follows from Assumption R5. \del{(\ref{ineq-inter}).}  This
completes the proof.

\section{Stochastic convolution in the cylindrical Gaussian case}

Assume now that $W(t)$, $t\geq 0$, is a cylindrical Wiener process
defined on some complete filtered probability space $(\Omega,\mathcal{F},(\mathcal{F}_t)_{t\geq 0},\mathbb{P})$. Let us
denote by $H$ the RKHS of that process, i.e. $H$ is equal to the
RKHS of $W(1)$.

\begin{theorem}\label{th-gaussian}
Under the above assumptions there exists a constant $\hat{C}_q(E)$
such that for any process $\xi$ described above the following
inequality holds
\begin{equation}
\label{ineq-2.05} \mathbb{E}  \int_0^T \vert
SC(\xi)(t)\vert_{D_A(\theta+\frac1p,q)} ^q \,dt \leq
\hat{C}_q(E)\mathbb{E} \int_0^T \Vert
\xi(t)\Vert_{R(H,D_A(\theta,q))} ^q \,dt , \; T\geq 0.
\end{equation}
\end{theorem}

The proof of Theorem \ref{th-gaussian} will be preceded by the
following useful result. \del{But first let us introduce a useful
notation.}

\begin{prop}\label{prop-interpolation_gaussian} Let us assume that $\theta
\in (0,1)$, $q \geq 1$ and $T>0$. Then there exists a constant
$K_T>0$ such that for each  bounded linear map $\vp:H\to E$ the
following inequality holds

\del{K_T^{-1} \Vert \vp \Vert_{R(H,D_A(\theta,q))}^q &\leq&
\nonumber\\\!\!\!\!\!\!\!\!\!\!\!\!\!\!\!\!&&\int_0^T
t^{(1-\vt)q}\Vert Ae^{-tA}\vp\Vert_{R(H,E)}^q \,\frac{dt}{t} \leq
K_T \Vert \vp \Vert_{R(H,D_A(\theta,q))}^q.}%

\begin{eqnarray}
\label{ineq-2.06} K_T^{-1} \Vert \vp
\Vert_{R(H,(E,D(A))_{\theta,q}}^q &\leq& \int_0^T t^{(1-\theta)
q}\Vert Ae^{-tA}\vp\Vert_{R(H,E)}^q \,\frac{dt}{t} \nonumber\\
&\leq& K_T \Vert \vp \Vert_{R(H,(E,D(A))_{\theta,q})}^q.
\end{eqnarray}

In particular, 
$\vp \in R(H,(D(A),E)_{\theta,q} )$ iff  (for some and/or all $T>0$)
the integral  $\int_0^T t^{(1-\theta) q}\Vert
Ae^{-tA}\vp\Vert_{R(H,E)}^q \,\frac{dt}{t}$ is finite.
\end{prop}
\begin{proof}[Proof of Proposition \ref{prop-interpolation_gaussian}.] Let
$\{e_k\}_k$ be an ONB of $H$ and $\{\beta_k\}_k$ a sequence of
i.i.d.\ Gaussian N(0,1) random variables. It is known, see e.g. \cite{Kw_Woy_1992} that there
exists a constant $C_p(E)$ such that for each linear operator $\vp
:H\to E$ the following inequality holds.
\begin{eqnarray}
\label{ineq-2.07} \hspace{0.5truecm} C_p(E)^{-1} \mathbb{E}\vert
\sum_j \beta_j \vp e_j\vert_E^p &\leq& \Vert \vp   \Vert_{R(H,E)}^p
\leq C_p(X)\mathbb{E}\vert \sum_j \beta_j \vp e_j\vert_E^p
\end{eqnarray}
\del{Denote $\theta=1-\vartheta$. Then} We have
\begin{eqnarray*}
&& \int_0^T t^{(1-\theta)q}\Vert Ae^{-tA}\vp\Vert_{R(H,E)}^q
\,\frac{dt}{t} \\&\leq& C_q(E) \int_0^T t^{(1-\theta)q}\mathbb{E}\vert
\sum_k \beta _k Ae^{-tA}\vp e_k\vert_{E}^q \,\frac{dt}{t}
\\
&=&C_q(E)\mathbb{E} \int_0^T t^{(1-\theta)q}\vert \sum_k \beta _k
Ae^{-tA}\vp e_k\vert_{E}^q \,\frac{dt}{t} \\&=& C_q(E) \mathbb{E}
\vert \sum_k \beta _k Ae^{-tA}\vp e_k\vert_{D_{A}(\theta ,q);T}^q
\\&\leq& C(T) C_q(E)\Vert \vp \Vert^q_{R(H,D_{A}(\theta ,q))}.
\end{eqnarray*}
Since $D_{A}(\vt ,q)=(E,D(A))_{\theta,q}$ with equivalent norms, this
proves the second inequality in (\ref{ineq-2.06}). The first inequality follows
the same lines.
\end{proof}

\begin{proof}[Proof of Theorem \ref{th-gaussian}.] From Proposition \ref{prop-interpolation_gaussian} we infer that the assumption (r5) in Theorem \ref{th-general} is satisfied. Since it is well known that the other assumptions are also satisfied, see e.g. \cite{Brz-03},
 the result follows from  Theorem \ref{th-general}.\end{proof}

\section{Proof of Theorem \ref{th-main}}

We only need to prove a version of Proposition \ref{prop-interpolation_gaussian} with $R(H,E)$ being replaced by $R(E):=L^p(S,\nu,E)$. We recall that here the measure space $(S,\CS,\nu)$ is fixed for the whole section.

\begin{prop}\label{prop-interpolation_poisson} Let us assume that $\theta
\in (0,1)$, $q \geq 1$ and $T>0$. Then there exists a constant
$K_T>0$ such that for each   $\vp\in L^p(S,\nu,E)=:R(E)$ the
following inequality holds

\del{K_T^{-1} \Vert \vp \Vert_{R(D_A(\theta,q))}^q &\leq&
\nonumber\\\!\!\!\!\!\!\!\!\!\!\!\!\!\!\!\!&&\int_0^T
t^{(1-\vt)q}\Vert Ae^{-tA}\vp\Vert_{R(E)}^q \,\frac{dt}{t} \leq
K_T \Vert \vp \Vert_{R(D_A(\theta,q))}^q.}%

\begin{eqnarray}
 K_T^{-1} \Vert \vp
\Vert_{R((E,D(A))_{\theta,q})}^q &\leq& \int_0^T t^{(1-\theta)
q}\Vert Ae^{-tA}\vp\Vert_{R(E)}^q \,\frac{dt}{t} \nonumber\\
&\leq& K_T \Vert \vp \Vert_{R((E,D(A))_{\theta,q})}^q.
\end{eqnarray}

In particular, 
$\vp \in R((D(A),E)_{\theta,q} )$ iff  (for some and/or all $T>0$)
the integral  $\int_0^T t^{(1-\theta) q}\Vert
Ae^{-tA}\vp\Vert_{R(E)}^q \,\frac{dt}{t}$ is finite.
\end{prop}
\begin{proof}[Proof of Proposition \ref{prop-interpolation_poisson}.] Follows by applying the Fubini Theorem.
\end{proof}

\appendix

\section{Martingale type $p$, $p\in [1,2]$, Banach spaces}
In this section we collect some basic information about the     martingale type $p$, $p\in [1,2]$, Banach spaces.

Assume also that $p\in [1,2]$ is fixed. A Banach space $E$ is of  martingale type $p$  iff
there exists a constant $L_{p}(E)>0$ such that for all
$X$-valued  finite martingale $\{M_{n}\}_{n=0}^N$  the
following inequality holds
\begin{equation} \sup_{n} \mathbb{E} | M_{n} | ^{p} \le  L_{p}(E)
\sum_{n=0}^N \mathbb{E}  | M_{n}-M_{n-1} | ^{p},
\label{eqn-2.1}\end{equation}
 where as  usually, we put  $M_{-1}=0$.

 Let us recall that a Banach space $X$ is of type $p$ iff there exists a constant $K_{p}(X)>0$  for any finite
sequence $\eps_1,\ldots,\eps_n \colon \Omega \to \{-1,1\}$ of
symmetric i.i.d. random variables and for any finite sequence
$x_1,\ldots,x_n$ of elements of $X$, the following inequality holds
\begin{equation}
 \mathbb{E} | \sum_{i=1}^n \eps_i x_i |^p \le
K_p(X) \sum_{i=1}^n | x_i |^p. \label{eqn-2.2}
\end{equation}

It is known, see e.g.\ \cite[Theorem 3.5.2]{0665.60005}, that  a Banach space $X$ is of type $p$ iff it is of Gaussian type $p$, i.e. there exists a constant ${\tilde K}_{p}(X)>0$  such that
for any finite
sequence $\xi_1,\ldots,\xi_n $ of  i.i.d.
$N(0,1)$ random  variables  and for any finite sequence
$x_1,\ldots,x_n$ of elements of $X$, the following inequality holds
         \begin{equation}
 \mathbb{E} | \sum_{i=1}^n \xi_i x_i |^p \le
\tilde{K}_p(X)  \sum_{i=1}^n | x_i |^p,
\label{eqn-2.3} \end{equation}

It is now well known, see e.g.\ Pisier \cite{Pis_1975} and  \cite{Pis_1986},  that $X$ is
of martingale type~$p$ iff it is $p$-smooth, i.e. there exists an equivalent  norm $|\cdot|$ on $X$ and there exist a constant $K>)$ such that
$\rho _{X}(t) \le  K t^{p}$ for all $t>0$, where  $\rho _{X}(t)$ is the
 modulus of smoothness of $(X,|\cdot|)$ defined by
$$ \rho _\mathcal{ X}(t) = \sup  \{ {1\over 2} ( | x+ty | + | x-ty | )-1
:  | x | , | y | =1\} .$$
In particular,  all spaces $L^{q}$ for $q\ge p$ and $q>1$, are of martingale type~$p$.

Let us also recall that  a  Banach  space $X$ it is
an  UMD  space  (i.e.  $X$ has   the unconditional  martingale
difference  property)  iff  for any $p \in (1,\infty)$ there exists a constant $\beta_p(X) > 0$ such that  for   any
$X$ -valued martingale difference $\{\xi _{j}\} $ ( i.e.:
$\sum^{n}_{j=1}\xi _{j}$  is a martingale), for any $\epsilon :
\mathbb{N} \to  \{-1,1\}$  and for any $n\in  \mathbb{N} $
\begin{equation}
\mathbb{E} |\sum^{n}_{j=1}  \epsilon_{j}\xi _{j} |^p \le  \beta_p
(X) \mathbb{E} |\sum^{n}_{j=1}  \xi_{j} |^p.
\label{eqn-2.6}\end{equation}

It is known,  see \cite{Burk_1986} and  references therein,  that  for  a
Banach space $X$ the following conditions are equivalent:
i) $X$ is an UMD space,
(ii) $X$ is $\zeta $  convex, \del{ i.e.,  there  is  a  biconvex
function $\zeta :X\times X \to   \mathbb{R} $
with the properties: $\zeta (0,0)>0, \zeta (x,y)\le  | x+y | $ for $
| x | , | y | =1$;} (iii) the Hilbert transform for $X$-valued functions is  bounded
in $L^{p}( \mathbb{R} ,X)$  for any (or some ) $ p>1$.

Finally, it is known, see e.g.\ \cite[Proposition 2.4]{Pis_1975}, that if a Banach space $X$ is both UMD and of type $p$, then $X$ is of martingale type $p$.

\section{Proof of inequality \ref{ineq-2.03}}\label{appB}
In this appendix we formulate and prove inequality \ref{ineq-2.03}. Our approach is a sense similar to the approach used in the Gaussian case by Neidhard \cite{Neidhardt-78} and Brze{\'z}niak \cite{Brz-95-M_type_2} or in the Poisson random measure in Madrekar and R\"udiger \cite{Man_Rud_2006}. In fact, our main result below can be seen a generalisation of Theorem 3.6 from \cite{Man_Rud_2006} to the case of martingale type $p$ Banach spaces.

\begin{Notation}
By $M_{S\times \mathbb{R}_+} ^{\bNN}$ we denote the family of all $\bNN$-valued measures on
$(S\times \mathbb{R}_+,\CS\otimes \mathcal{B}_{\mathbb{R}_+})$ and
$\CM_{S\times \mathbb{R}_+} ^ {\bNN}$ is the $\sigma$-field on $M_{S\times \RR_+} ^{ \bNN}$ generated by functions
$i_B:M\ni\mu \mapsto \mu(B)\in \bNN$, $B\in \CS\otimes \mathcal{B}_{\mathbb{R}_+}$.
\end{Notation}

Let us assume that  $(S,\CS)$ is  a measurable space, $\nu\in M ^ +_S$ is a non-negative measure on $(S,\CS)$ and  $\mathfrak{P}=(\Omega,\CF,(\mathcal{F}_t)_{t\geq 0},\PP)$ is a filtered
probability space. We also assume that  $\eta$ is time homogeneous Poisson random measure  over $\mathfrak{P}$,  with the intensity measure $\nu$, i.e.   $\eta: (\Omega,\CF)\to (M_S ^ \bNN,\CM ^ \bNN_{S\times \mathbb{R}_+}) $
is a measurable function  satisfying  the following conditions
\begin{trivlist}
\item[(i)] for each $B\in  \CS \otimes \mathcal{B}_{\mathbb{R}_+} $,
 $\eta(B):=i_B\circ \eta : \Omega\to \bar{\mathbb{N}} $ is a Poisson random variable with parameter\footnote{If  $\EE \eta(B) = \infty$, then obviously $\eta(B)=\infty$ a.s..} $\EE\eta(B)$;
\item[(ii)] $\eta$ is independently scattered, i.e. if the sets
$ B_j \in   \CS\otimes \mathcal{B}_{\mathbb{R}_+}$, $j=1,\cdots, n$ are pair-wise disjoint,   then the random
variables $\eta(B_j)$, $j=1,\cdots,n $ are pair-wise independent;
\item[(iii)] for all $B\in  \CS $ and $I\in \mathcal{B}_{\mathbb{R}_+}$, $\mathbb{E}\big[\eta (I\times B)\big]=\lambda(I)\nu(B)$, where $\lambda$ is the Lebesgue measure;
\item[(iv)] for each $U\in \CS$, the $\bar{\mathbb{N}}$-valued processes $(N(t,U))_{t>0}$  defined by
$$N(t,U):= \eta((0,t]\times U), \;\; t>0$$
is $(\mathcal{F}_t)_{t\geq 0}$-adapted and its increments are independent of the past, i.e. if $t>s\geq 0$, then
$N(t,U)-N(s,U)=\eta((s,t]\times U)$ is independent of
$\mathcal{F}_s$.
\end{trivlist}
By  $\tilde{\eta}$  we will denote the compensated Poisson random measure, i.e. a function defined by $\tilde\eta(B)=\eta(B)-\mathbb{E}(\eta(B))$, whenever the difference makes sense.

\del{Moreover, let us recall the predictable $\sigma$ field which is the $\sigma$ field on $\Omega\times \RR_+$ generated by
all sets $A\in \CF\hat \times \CB( \RR_+)$, where $A$ is of the form $A=F\times (r,t]$,
with $F\in\CF_r$ and $r,t\in\RR_+$.
If $\xi:\Omega\times \RR_+\to S$ is $\CP$ measurable, then $\xi$ is called a predictable process (with respect to the filtration $(\mathcal{F}_t)_{t\geq 0}$).}

\begin{lem}
\label{lem-cw5}
Let $p\in(1,2]$ and assume that $E$ is a Banach space of martingale type $p$.   If  a finitely-valued function  $f$ belongs to  $ L^p(\Omega\times S,\mathcal{F}_{a}\otimes \CS ;\mathbb{P}\otimes \nu;E)$ for some $a\in \mathbb{R}_+$, then
for any $b >a$,
\begin{equation}
\mathbb{E}\vert \int_S f (x)\tilde{\eta}(dx,(a,b]) \vert_E^p
\leq
2^{2-p}L_p(E)  (b-a) \mathbb{E}
\int_S \vert f(x)\vert_E^p \, \nu(dx) 
\end{equation}
\end{lem}
Since the space of finitely-valued functions  is dense in   $ L^p(\Omega\times S,\mathcal{F}_{a}\otimes \CS ;\mathbb{P}\otimes \nu;E)$, see e.g. Lemma 1.2.14 in \cite{Chalk-06}.

\begin{coro}Under the assumptions of Lemma \ref{lem-cw5} there exists a unique bounded linear operator
$$
\tilde{I}_{(a,b)}:  L^p(\Omega\times S,\mathcal{F}_{a}\otimes \CS ;\mathbb{P}\otimes \nu;E) \to L^p(\Omega,\mathcal{F},E)
$$
such that for a finitely-valued function $f$,
we have
$$
\tilde{I}_{(a,b)}(f)=\int_S f (x)\tilde{\eta}(dx,(a,b]).
$$
In particular, for every $f\in L^p(\Omega\times S,\mathcal{F}_{a}\otimes \CS ;\mathbb{P}\otimes \nu;E)$,
\begin{equation}
\mathbb{E}\vert \tilde{I}_{(a,b)}(f) \vert_E^p
\leq
2^{2-p}L_p(E)  (b-a) \mathbb{E}
\int_S \vert \xi(x)\vert_E^p \, \nu(dx).
\end{equation}
In what follows, unless we in danger of ambiguity, for every
$L^p(\Omega\times S,\mathcal{F}_{a}\otimes \CS ;\mathbb{P}\otimes \nu;E)$
we will write $\int_S \xi (x)\tilde{\eta}(dx,(a,b])$ instead of $\tilde{I}_{(a,b)}(f)$.
\end{coro}

Let $X$  be any Banach	space. Later on we will take $X$ to  be one of the
spaces $E$, $R(H,E)$ or  $L^p(S,\nu,E)$.
For $a<b\in [0,\infty]$ let
$\mathcal{N}(a,b;X)$ be the space of (equivalence classes of)
predictable functions $\xi :(a,b]\times\Omega\to X$.

For $ q\in (1,\infty ) $ we  set
\begin{eqnarray}\;\;
\mathcal{N}^q(a,b;X)&=&
\left\{
\xi \in \mathcal{N}(a,b;X): \;
 \int_a^b\vert\xi (t)\vert^q\,dt<\infty \mbox{ a.s. }
\right\},
\label{def:Nq}
\\
\;\; \;\;
\mathcal{M}^q(a,b;X)&=&
\left\{\xi \in \mathcal{N}(a,b;X):\mathbb{E}\int_
a^b\vert\xi (t)\vert^q\,dt<\infty\right\}.
\label{def:Mq}
\end{eqnarray}
Let $\mathcal{N}_{\rm step}(a,b;X)$ be the space of all $
\xi \in \mathcal{N}(a,b;X)$ for which  there exists a partition
$a=t_0<t_1<\cdots <t_n<b$ such that for $k\in\{1,\cdots,n\}$, for $t\in (t_{k-1},t_{k}]$,
$\xi (t)=\xi (t_k)$ is $\mathcal{F}_{t_{k-1}}$-measurable and $\xi(t)=0$ for $t\in (t_n,b)$. We put $\mathcal{M}_{\mbox{\small step}}^q
=\mathcal{M}^q\cap \mathcal{N}_{\mbox{\small step}}$.
Note that $\mathcal{M}^q(a,b;X)$ is a closed subspace of $L^q([a,b)\times
\Omega;X)\cong L^q([a,b);L^q(\Omega;X))$.

In what follows we put $a=0$ and $b=\infty$.
For  $\xi\in \mathcal{M}_{\mbox{\small step}}^p(0,\infty;L^p(S,\nu;E))$ we set
\begin{equation}
 {\tilde I}(\xi )=\sum_{j=1}^n\int_S \xi (t_j,x)\tilde{\eta}(dx , (t_{j-1},t_j]).
\label{eqn-int-def}
\end{equation}

 \del{Since  for $\omega\in\Omega$,
$\xi (t_k,\omega )\in \mathcal{L}(E,X)$
and  $w(t_{k+1},\omega)-w(t_k,\omega)\in E$,} Obviously,   ${\tilde I}(\xi )$
is a $\mathcal{F}$-measurable map from $\Omega$ with values in $E$.

We have the
following  auxiliary  results.
\begin{lem}
\label{lem-cw4}
Let $p\in(1,2]$ and assume that $E$ is a Banach space of martingale type $p$.  Then for any  $\xi\in \mathcal{M}_{\mbox{\small step}}^p(0,\infty;L^p(S,\nu;E))$,
${\tilde I}(\xi )\in L^p(\Omega ,E)$,
$\mathbb{E}{\tilde I}(\xi
)=0$ and
\begin{equation}
\mathbb{E}\vert {\tilde I}(\xi )\vert^p
\leq
L^2_p(E) 2^{2-p} \int_0^\infty \mathbb{E}
\int_S \vert \xi(t,x)\vert_E^p \nu(dx) \,dt
\label{I1}
\end{equation}
\end{lem}
%
\begin{lem}
\label{lem-cw6}
Suppose that $\xi \sim \text{ Poiss }(\lambda)$,
where $\lambda >0$. Then, for all $p\in [1,2]$,
$$\mathbb{E}|\xi-\lambda|^p \leq 2^{2-p}\lambda.$$
\end{lem}

\begin{remark}
One can  easily calculate that
$$\mathbb{E}(|\xi-\lambda|)=2\lambda\, e^{-\lambda},\;\text{ if }
\lambda \in (0,1).$$
\end{remark}
\del{ Let us recall the predictable $\sigma$ field which is the $\sigma$ field on $\Omega\times \RR_+$ generated by
all sets $A\in \CF\hat \times \CB( \RR_+)$, where $A$ is of the form $A=F\times (r,t]$,
with $F\in\CF_r$ and $r,t\in\RR_+$.
If $\xi:\Omega\times \RR_+\to S$ is $\CP$ measurable, $\xi$ is called predictable. Moreover, let  $\mathcal{M}^p(0,\infty, L^p(S,\nu;E))$ be the space of all
predictable processes $\xi:[0,\infty)\to  L^p(S,\nu;E)$, such that
\DEQSZ\label{norm-cm}
 \lqq{\phantom{mmmmmmmmmm}|x|_{\CM} := }
\\
\nonumber &&\lk( \int_0 ^ \infty \EE\int_S |\xi(t,x)|_E ^ p\nu(dx)\: dt\rk) ^ \frac 1p<\infty , \quad x\in \mathcal{M}^p(0,\infty;L^p(S,\nu;E)).
\EEQSZ
\begin{remark}
The space $\mathcal{M}^p(0,\infty, L^p(S,\nu;E))$ denotes the completion of $\CM_{\mbox{\small step}}^p(0,\infty;L^p(S,\nu;E))$ with respect to the norm \eqref{norm-cm}.
\end{remark}
}

\begin{theorem}\label{th-ineq} Assume that $p\in (1,2]$ and $E$ is a martingale type $p$ Banach space. Then there exists
there exists a unique bounded linear operator
$$
\tilde{I}:  \mathcal{M}^p(0,\infty, L^p(S,\nu;E)) \to L^p(\Omega,\mathcal{F},E)
$$
such that for
$\xi \in \mathcal{M}_{\rm step}^p(0,\infty, L^p(S,\nu;E)) $ we have  $I(\xi)=\tilde{I}(\xi)$.  In particular, for every $\xi \in \mathcal{M}^p(0,\infty, L^p(S,\nu;E))$,
\begin{equation}
\mathbb{E}\vert {I}(\xi) \vert_E^p
\leq
2^{2-p}L_p^2(E)   \mathbb{E}
\int_0^\infty \int_S \vert \xi(t,x)\vert_E^p \, \nu(dx)\, dt.  \label{I4}
\end{equation}
\end{theorem}

\begin{proof}[Proof of Theorem \ref{th-ineq}.] Follows from Lemma \ref{lem-cw4} and the density of $\mathcal{M}_{\rm step}^p(0,\infty, L^p(S,\nu;E))$ in the space $\mathcal{M}^p(0,\infty, L^p(S,\nu;E))$.
\end{proof}

In a natural way we can define spaces $\mathcal{M}_{\rm loc}^p(0,\infty, L^p(S,\nu;E))$ and $\mathcal{M}^p(0,T, L^p(S,\nu;E))$, where $T>0$. Then for any $\xi \in \mathcal{M}_{\rm loc}^p(0,\infty, L^p(S,\nu;E))$ we can in a
standard way define the  integral $\int_0^t \int_S\xi(r,x)\tilde{\eta}(dx,dr)$, $t\geq 0$, as the c\'adl\'ag modification of the process
\begin{eqnarray}
I(1_{[0,t]}\xi),\; t\geq 0,
\end{eqnarray}
 where  $[1_{[0,t]}\xi](r,x;\omega):= 1_{[0,t]}(r)\xi(r,x,\omega)$, $t\geq 0$, $r\in\RR_+ $, $x\in S$  and $\omega\in\Omega$.
To show that this c\'adl\'ag modification exists we argue as follows. First of all we can assume that  $\mathcal{M}^p(0,T, L^p(S,\nu;E))$, for some $T>0$. Let  $\{\xi_n\}_{n\in\NN}$
be an $\mathcal{M}_{\rm step}^p(0,T, L^p(S,\nu;E))$-valued  sequence that is  convergent in $\mathcal{M}^p(0,T, L^p(S,\nu;E))$ to $\xi$. Hence, the sequence  $\{\xi_n,n\in\NN\}$ is uniformly integrable  and so it follows that the condition (a) in  Remark 3.8.7 from \cite{ethierkurtz} is satisfied.
\del{Observe, one has to take as underlying probability space $(\Omega\times[0,T])$ equipped by the
product measure $\PP_T$ defined by
$$\PP_T ((A,B))= \int_A \: dt\; \PP(B),
$$
and to consider the uniform integrability condition with respect to $\PP_T$ on $(\Omega\times[0,T])$.}
Similarly, the compact containment condition, i.e. the condition (a) in  Theorem 3.7.2  from \cite{ethierkurtz}, holds true in view of the Prohorov Theorem,  since for any $t\ge 0$ the laws of the sequence
$\{ I(1_{[0,t]}\xi_n) ,n\in\NN\}$ are tight in the set of all probability measures over $E$, compare also with \cite{dett}.

Similarly, for a stopping time $\tau$ we can define and process $\xi \in \mathcal{M}_{\rm loc}^p(0,\infty, L^p(S,\nu;E))$ and  the  integral
\begin{eqnarray}
\int_0^\tau \int_S\xi(r,x)\tilde{\eta}(dx,dr):=I(1_{[0,\tau]}\xi),
\end{eqnarray}
provided $1_{[0,\tau]}\xi \in \mathcal{M}^p(0,\infty, L^p(S,\nu;E))$. Theorem \ref{th-ineq} implies that in this case the following inequality holds.
\begin{equation}
\mathbb{E}\vert \int_0^\tau \int_S\xi(r,x)\tilde{\eta}(dx,dr) \vert_E^p
\leq
C_p   \mathbb{E}
\int_0^\tau \int_S \vert \xi(r,x)\vert_E^p \, \nu(dx)\, dr.  \label{I4bis}
\end{equation}
with some constant $C_p>0$ independent of $\xi$.

\begin{proof}[Proof of Lemma \ref{lem-cw4}.] Let us observe that the sequence $(M_k)_{k=1}^n$ defined by $M_k=\sum_{j=1}^k\int_S \xi (t_j,x)\tilde{\eta}(dx,[t_{j-1},t_j))$ is an $E$-valued martingale (with respect to the filtration $(\mathcal{F}_{t_k})_{k=1}^n$). Therefore, by the martingale type $p$ property of the space $E$ and Lemma \ref{lem-cw5} we have the following sequence of inequalities
\begin{eqnarray}\nonumber
\mathbb{E}\vert(\tilde{I}(\xi))\vert^p_E &=& \mathbb{E}\vert M_n\vert^p_E\leq L_p(E) \sum_{k=1}^n \mathbb{E} \vert \int_S \xi (t_k,x)\tilde{\eta}(dx,[t_{k-1},t_k]) \vert_E^p\\
&\leq& L^2_p(E)2^{2-p} \sum_{k=1}^n (t_k-t_{k-1})\: \mathbb{E}  \int_S \vert \xi(t_k,x)\vert_E^p\,\nu(dx)\\
\nonumber
&=& L^2_p(E)2^{2-p} \int_0^\infty \mathbb{E}
\int_S \vert \xi(t,x)\vert_E^p \, \nu(dx) \,dt.
\end{eqnarray}
This concludes the proof.
\end{proof}

\begin{proof}[Proof of Lemma \ref{lem-cw5}.] Put $I=(a,b]$. We may
suppose that $f=\sum_{i}f_i1_{A_i\times B_i}$  with $f_i\in E$, $A_i\in \mathcal{F}_{a}$ and $B_i\in\CS$,
 the finite family of sets $({A_i\times B_i})$ being pair-wise disjoint and $\nu(B_i)<\infty$.  Let us notice that
$$
  \int_S f (x)\tilde{\eta}(dx,I)=\sum_{i}1_{A_i} \tilde \eta({B_i\times I }) f_i.$$

 Since the random variables $\tilde{\eta}({B_i \times I})$ are independent from the $\sigma$-field $\mathcal{F}_a$, the random variables
  $1_{A_i}\tilde{\eta}({B_i\times I })$ conditioned on $\mathcal{F}_a$ are independent and so by the martingale type $p$ property of the space $E$ and Lemma \ref{lem-cw5} we infer that
\begin{eqnarray*}
 \mathbb{E}\vert \int_S \xi (x)\tilde{\eta}(dx,I)\vert_E^p&=&\mathbb{E}\big[\mathbb{E} \big(\vert \sum_{i}1_{A_i} \tilde \eta({B_i\times I}) f_i \vert_E^p |\mathcal{F}_a \big)\big]
 \\
 &\leq&  \mathbb{E}\big[ L_p(E) \sum_{i} \vert f_i 1_{A_i}\vert_E^p \mathbb{E}|\tilde{\eta}( B_i\times I)|^p \big]
 \\&\leq & L_p(E) \mathbb{E}\big[  \sum_{i} \vert f_i \vert_E^p 1_{A_i} 2^{2-p}\lambda(I)\nu(B_i) \big]
 \\
 &=& 2^{2-p}L_p(E)\sum_{i}\vert f_i\vert ^p_E
 \nu(B_i) \lambda(I)  \mathbb{P}(A_i)\\
 &=& 2^{2-p}L_p(E) \lambda(I) \int_{\Omega\times S} \vert \sum_{i} f_i 1_{A_i\times B_i}\vert^p\, d(\mathbb{P}\otimes \nu)\\
 &=&\tilde L_p(E)(b-a) \mathbb{E}
\int_S \vert f(x)\vert_E^p \, \nu(dx).
 \end{eqnarray*}
 The proof is complete.
\end{proof}
\begin{proof}[Proof of Lemma \ref{lem-cw6}.] The case $p=2$ is well known. Since $\xi \geq 0$ and
$\mathbb{E}(\xi)=\lambda$, the case $p=1$ follows by the triangle
inequality. The case $p\in (1,2)$ follows then by applying the
H\"older inequality. Indeed, with $\alpha=2(p-1)$ and $\beta=2-p$ we
have the following sequence of inequalities, where
$\eta:=|\xi-\lambda|$.
\begin{eqnarray*}
\mathbb{E}(\eta^p) &=&\mathbb{E}(\eta^\alpha \eta^\beta)  \leq
[\mathbb{E}((\eta^\alpha)^{2/\alpha})]^{\alpha/2}
[\mathbb{E}((\eta^\beta)^{1/\beta})]^{\beta}\\&=&
[\mathbb{E}(\eta^2)]^{\alpha/2}[\mathbb{E}(\eta)]^{\beta} \leq
(\lambda)^{\alpha/2}(2\lambda)^{\beta}=2^{2-p}\lambda.
\end{eqnarray*}
\end{proof}

We conclude with a result corresponding to inequality \eqref{ineq-2.04}.

\begin{coro}\label{cor_maximal} Assume that $1<q\leq p<2$ and $E$ is a martingale type $p$ Banach space. Then there exists
there exists a constant $C>0$ such that for any process $\xi \in \mathcal{M}^p_{\rm loc}(0,\infty, L^p(S,\nu;E)) \to L^p(\Omega,\mathcal{F},E)$, and any $T>0$,
\begin{equation}
\label{ineq-2.04bis} 
\mathbb{E} \vert \sup_{t\in [0,T]} \int_0^t \int_S \xi(r,x)
\tilde{\eta}(dx,dr)\vert ^q \leq C\del{C_q(E)}\mathbb{E}
\big(\int_0^T\int_S\vert \xi(r,x)\vert^p\,
\nu(dx)\,dr\big)^{q/p}.
\end{equation}
\end{coro}

The proof of the above result will be based on Proposition IV.4.7 from the monograph \ref{prop-RY} by Revuz and Yor  which we recall here for the convenience of the reader.
\begin{prop}
\label{prop-RY} Suppose that a positive, adapted right-continuous
process $Z$ is dominated by an increasing process $A$, with $A_0$,
i.e. there exists a constant $C>0$ such that for every bounded stopping time $\tau$, $\mathbb{E}Z_\tau \leq
C \mathbb{E}A_\tau$. Then for  any $k \in (0,1)$,
$$
\mathbb{E} \sup_{0\le t < \infty} Z_t^k \leq C^k \frac{2-k}{1-k}
\mathbb{E}A_\infty^k.
$$
\end{prop}

\begin{proof}[Proof of Corollary \ref{cor_maximal}] Let now fix $q \in (1,p)$. Put  $k=q/p$. We will apply Proposition \ref{prop-RY} to
the processes $Z_t=\vert \int_0^t \int_S \xi(r,x)
\tilde{\eta}(dx,dr) \vert_E^p$ and $A_t=
\int_0^t \int_S \vert \xi(r,x)\vert_E^p \, \nu(dx)\, dr
$, $t\in [0,T]$. Let us notice that in view of  inequality \eqref{I4bis}, the process  $Z$ is dominated by the process $A$. Since
$Z$ is right continuous, $\sup_{0\le t
\leq T} Z_t^k= \sup_{0\le t \leq T}\vert \int_0^t \int_S \xi(r,x)
\tilde{\eta}(dx,dr) \vert_E^q $ and
$A_\infty^k =\left( \int_0^T \int_S \vert \xi(r,x)\vert_E^p \, \nu(dx)\, dr\right)^{q/p}$,
we get inequality \eqref{ineq-2.04bis}.
This completes the proof of Corollary  \ref{cor_maximal}. \end{proof}

\def\cprime{$'$} \def\cprime{$'$} \def\cprime{$'$} \def\cprime{$'$}

\vskip 2truecm



\end{document}